\documentclass[12pt]{amsart}
\usepackage{graphicx}
\usepackage{pb-diagram}
\newcommand\lra{\longrightarrow}
\newcommand\PP[1][3]{\mathbb P^{#1}}
\newcommand{\ZZ}{\mathbb{Z}}
\newcommand{\QQ}{\mathbb{Q}}
\newcommand{\RR}{\mathbb{R}}
\newcommand{\CC}{\mathbb{C}}
\newcommand{\FF}{\mathbb{F}}
\def\arraystretch{1.2}
\theoremstyle{remark}

 \def\dom#1,#2,#3,#4,#5,#6{#1\le x\le#2,\;#3\le y\le#4,%
\; #5\le z\le#6}

\begin{document}
\title{Periods of double octic Calabi--Yau manifolds}
\author{S. Cynk}
\address{S\l awomir Cynk:
Institute of Mathematics, Jagiellonian University
\L ojasiewicza 6,
30-348 Krak\'ow,
POLAND}
\address{
Institute of Mathematics of the
Polish Academy of Sciences, ul. \'Sniadeckich 8, 00-956 Warszawa, POLAND
}
\email{slawomir.cynk@uj.edu.pl}
\author{D. van Straten}
\address{Duco van Straten: Fachbereich 08,
AG Algebraische Geometrie,
Johannes Gu\-ten\-berg-Universit\"at,
D-55099 Mainz, GERMANY}
\email{straten@mathematik.uni-mainz.de}

\begin{abstract}
  We compute numerical approximations of the period integrals for eleven rigid
  double octic Calabi--Yau threefolds and compare them with the
  periods of corresponding weight four cusp forms and find, as to be expected,
commensurabilities. These give information on character of the correspondences
of these varieties with the associated Kuga-Sato modular threefolds.
\end{abstract}
\maketitle

\section*{Introduction}
Let $X$ be a rigid Calabi--Yau threefold and $\omega\in H^{3,0}(X)$ a
regular $3$--form on $X$. For a $3$-cycle $\gamma \in H_3(X,\ZZ)$ on $X$ 
we can form
the {\em period integral}
\[ \int_{\gamma} \omega \in \CC .\]
The set of these period integrals form a lattice
\[\Lambda:=\left\{\int_{\gamma}\omega:\;\gamma\in
  H_{3}(X,\ZZ)\right\} \subset \CC\]
and hence determine an elliptic curve 
\[ \CC/\Lambda={H^{3,0}(X)}^{*}/H_{3}(X,\ZZ)=:J^2(X)\]  
which is just an example of the {\em intermediate Jacobian} of
Griffiths \cite{Griff}.\\

It is now known that any rigid Calabi--Yau threefold defined over
$\QQ$ is {\em modular} \cite{GY},\cite{Die}, in the sense 
that one has an equality of $L$-functions: 
\[ L(H^3(X),s)=L(f,s) .\]
Here $f \in S_4(\Gamma_0(N))$ is a weight four cusp form for some level $N$.
It is known more generally that a cusp-form $f \in S_{k+2}(\Gamma_0(N))$ can be
interpreted as a $k$-form on the associated Kuga-Sato variety, which
is (a desingularisation of) the $k$-fold fibre product of the
universal elliptic curve over the modular curve $X_0(N)$,
\cite{Del}. 
So in our case one expects the   
equality of $L$-functions to come from a correspondence between the
rigid  
Calabi-Yau and the Kuga-Sato variety $Y$, which resolves the fibre product
$ E \times_C E$ of the universal elliptic curve over the modular curve
$C:=X_0(N)$.
Now the periods
\[ \int_0^{i\infty} f(\tau)\tau^k d\tau\]
of the modular form determine the periods of the Kuga-Sato variety $Y$
and the correspondence between $X$ and  $Y$ would  imply that the 
period-lattice of $X$ is commensurable to the lattice derived from the 
modular form.\\

In this note we shall compute numerical approximations of period
integrals for certain rigid double octic Calabi--Yau threefolds,
i.e. Calabi--Yau threefolds constructed as a resolution of a double
cover of projective space $\mathbb P^{3}$, branched along a surface
of degree eight.\\

More specifically, we will look at the eleven arrangements of
eight planes defined by linear forms with rational coefficients, 
described in the PhD thesis of {\sc C. Meyer} \cite{Meyer}. These 
arrangements define eleven non--isomorphic rigid Calabi--Yau threefolds. 
He also determined the weight four cusp forms $f$ for these eleven rigid 
double octics using the counting of points in $\FF_p$ for small primes $p$.\\

Each arrangement of real planes defines a partition of the real projective space $\mathbb{P}^3(\RR)$ into {\em polyhedral cells} and using these cells one can construct certain {\em polyhedral $3$-cycles} on the desingularisation of the double octic. Using the explicit equations for the planes of the arrangement,
one can write the period integral as an explicit sum of multiple integrals, 
which can be integrated numerically.\\

It turned out to be difficult to identify a complete basis of $H_3(X,\ZZ)$ in 
terms of polyhedral cycles. But any two non--proportional periods of a rigid 
Calabi--Yau threefolds define a subgroup of finite index of $\Lambda$ and 
hence an elliptic curve {\em isogeneous} to the intermediate jacobian $J^2(X)$. 
From such a numerical lattice one can compute the lattice constants $g_2$ and $g_3$, and hence a Weierstrass equation and $j$--invariant of the curve defined by lattice  spanned by the polyhedral $3$--cycles.\\ 

We expect that a more refined topological analysis of the above situation 
will lead to more precise information on the nature of the correspondences 
between these varieties.\\ 

\section{Double octic Calabi--Yau threefolds}
\label{sec:doct}
By a {\em double octic} we understand a variety $X$ given as a double cover
\[ \pi: X \lra \PP\]
of $\PP$, ramified over a surface $D \subset \PP$ of degree eight.
Such a double octic $X$ can be given by an equation in weighted
projective space $\PP[](4,1,1,1,1)$ of the form
\[ u^2=F(x,y,z,t) ,\]
where the polynomial $F$ defines the ramification divisor $D$. If the surface 
$D$ is smooth, then $X$ is a smooth Calabi-Yau threefold, but we will be 
dealing here with the case that $D$ is a union of eight planes, so the 
polynomial actors as into a product of linear forms:
\[ F=L_1 L_2 L_3 L_4 L_5 L_6 L_7 L_8 .\]
The associated double octic $X$ then is singular along the lines of
intersection of the eight planes $D_i:=\{L_i=0\}$. In case these
planes have the property that
\begin{center}
\begin{itemize}
\item[] no six intersect in a point,
\item[] no four intersect along a line
\end{itemize}
\end{center}
one can construct a Calabi-Yau desingularisation of $X$.
To do so, one first constructs a sequence of blow--ups with  smooth 
centers 
\[ T: \widetilde{\PP} \lra \PP\]  
and a divisor $\widetilde D$ in $\widetilde \PP$ such that
\begin{itemize}
\item[] $\widetilde D$ is {\em non--singular} (in particular reduced),
\item[] $\widetilde D$ is {\em even} as an element of the Picard group
  $\operatorname{Pic}(\widetilde{\PP})$,
\end{itemize}
by blowing--up the singularities of $D$ in the following order:
\begin{center}
\begin{enumerate} 
\item fivefold points,
\item triple lines,
\item fourfold points,
\item double lines.
\end{enumerate}
\end{center}
In the first two cases we replace the branch divisor by its reduced inverse
image: the strict transform plus the exceptional divisor. In the last
two cases we replace the branch divisor by its strict transform.

The double cover 
\[\widetilde \pi:\widetilde X  \lra \widetilde{\PP}\]
of $\widetilde \PP$ branched along $\widetilde D$ is now a smooth Calabi--Yau
manifold, which we will call the {\em double octic Calabi-Yau threefold}
of the arrangement.\\

We will also need to consider a particular {\em partial resolution} $\hat X$
of $X$, obtained as double cover of a space $\hat \PP$ obtained by
performing only the blow--ups in fivefold point, triple lines and
double 
curves, so leaving out step $(3)$ in the above procedure. Note that 
there are two types of fourfold points. The fourfold points on triple
lines get removed in step $(2)$, but the fourfold points that appear at
the intersection of four generic planes produce an {\em ordinary double point}
if we blow up consecutively the six curves of intersection of the
strict transforms of these planes. 
After the blow-up of the first double line the strict transforms of
the remaining two planes (not containing this line) intersect along
a sum of two intersecting lines. So we have to blow-up four lines and a
cross, the latter producing a node on the threefold $\hat \PP$.  
The space doubly covering $\hat \PP$ and ramified over $\hat D$ is a variety
$\hat X$ with twice as many nodes.\\

To summarise the situation, one can consider the following diagram:
\[
  \begin{diagram}\dgARROWLENGTH=1.2cm
    \node{\tilde X} \arrow{e,t}{\sigma} \arrow{s,l}{\tilde\pi}
    \node{\hat X} \arrow{e,t}{\rho}\arrow{s,l}{\hat\pi}\node {X}
    \arrow{s,r}{\pi}\\
    \node{\widetilde \PP}\arrow{e,t}{S} \node{\widehat \PP}
    \arrow{e,t}{R} \node{\PP}
  \end{diagram}
\]
The vertical maps are two-fold covers, the map $\sigma: \widetilde X
\lra \hat X$ is a small resolution of the nodes of $\hat X$ and $\rho:
\hat X \lra X$ is a 
partial resolution of double octic variety $X$. The composition $R S$ is
the map $T:\widetilde \PP \lra \PP$ we started with, and $\tau:=\rho
\sigma: \widetilde X \lra X$ is a resolution of singularities.
In fact the resolutions $\tilde X$ and $\hat X$ depend on the choice of the
order of blow-up of lines, in the above diagram we 
choose the same order of lines for both resolutions. 
\\

We will be concerned with $11$ special arrangements that were studied by
{\sc C. Meyer} \cite{Meyer}. Their resolution $\widetilde X$ of the associated
double octic lead to $11$ different rigid Calabi-Yau varieties. For the
convenience of the reader, we list here the arrangement numbers,
second Betti-number and the equations from \cite{Meyer}. 
\[\scriptsize
\begin{array}{|c|c|c|r|}
\hline\rule[-0.5mm]{0mm}{4mm}
\textup{Number}&b_2(\tilde X)&\textup{Equation}&\lambda\ \\
\hline
1&70&xyzt(x + y)(y + z)(z + t)(t + x)&-1\\
\hline
3&62&xyzt(x + y)(y + z)(y -t)(x - y -z + t)&1\\
\hline
19&54&xyzt(x + y)(y + z)(x - z - t)(x + y + z - t)&2\\
\hline
32&50&xyzt(x + y)(y + z)(x - y - z - t)(x + y - z + t)&-1\\
\hline
69&50&xyzt(x + y)(x - y + z)(x - y - t)(x + y - z - t)&-1\\
\hline
93&46&xyzt(x + y)(x - y + z)(y - z - t)(x + z - t)&2\\
\hline
238&44&xyzt(x + y + z - t)(x + y -z + t)(x - y + z + t)(-x + y + z + t)&1\\
\hline
239&40&xyzt(x + y + z)(x + y + t)(x + z + t)(y + z + t)&1\\
\hline
240&40&xyzt(x + y + z)(x + y - z + t)(x - y + z + t)(x - y - z - t)&-2\\
\hline
241&40&xyzt(x + y + z + t)(x + y - z - t)(y - z + t)(x + z - t)&1\\
\hline
245&38&xyzt(x + y + z)(y + z + t)(x - y - t)(x - y + z + t)&-2\\
\hline
\end{array}
\] 
The meaning of the number $\lambda$ in the last column will be explained later.\\

\section{3--cycles on a double octic}

In the above eleven examples of rigid double octic Calabi--Yau
threefolds defined over $\mathbb Q$ are given. In all these
examples the eight planes are given by equations with integral 
coefficients.\\

In general, an arrangement defined by real planes gives a decomposition of
$\PP(\mathbb R)$ into a finite number of {\em polyhedral cells}. By combining
these cells one can construct certain {\em polyhedral cycles} on the
smooth model $\widetilde X$.
To explain this, let us fix one of these cells $C$ and consider its double covering ${\bf C}$, that is, its preimage under the $2$--fold covering map 
$\pi: X \lra \PP$.
Then ${\bf C}$ is a $3$--cycle in $X$ and determines an element in $H_3(X,\ZZ)$,
up to a sign determined by a choice of an orientation.\\

{\bf Question:} Do the $3$--cycles ${\bf C}$ generate $H_3(X,\ZZ)$?\\

However, ${\bf C}$ will in general not be a $3$--cycle on the desingularisation $\widetilde X$, as the canonical map 
\[ \tau_*: H_3(\widetilde X,\ZZ) \lra H_3(X,\ZZ)\]
will not be surjective in general. To see this geometrically, we follow the 
fate of the cycle ${\bf C}$ under the blow--up maps and see that its gets 
transformed into a chain on $\widetilde X$ that we will still denote by ${\bf C}$. Its boundary $\partial {\bf C}$, as a
chain on $\widetilde X$, is a sum of $2$--cycles contained in the exceptional loci
\[\partial {\bf C}=\bigcup_{i}\Gamma_{i}.\] 
Now observe that ${\bf C}$ is {\em anti--symmetric} with respect to the 
covering map
$\pi$, and as a consequence, these cycles $\Gamma_i$ are
{anti--symmetric} as well.
On the other hand the exceptional divisor corresponding to a fivefold
point or a triple line is fixed by the involution, while the
exceptional locus corresponding to a double line is a blow-up of a
conic bundle, so in all the three cases second cohomology group is
symmetric. 
Hence  each cycle $\Gamma_{i}$ contained in
the  
exceptional divisor corresponding to a double line, triple line or a fivefold 
point there is a boundary,
i.e. there is a $3$--chain ${\bf C}_{i}$ such that $\delta {\bf C}_{i}=\Gamma_{i}$. 
Hence, if we subtract from the chain ${\bf C}$ the chains ${\bf C}_{i}$ we get
a chain with boundary contained in the exceptional divisors
corresponding to the fourfold points. So we see that ${\bf C}$ can be
lifted to a cycle $\bf \widetilde{C}$ on the partial resolution $\hat X$ and we have shown:\\

{\bf Proposition:} The map
\[ \rho_*: H_3(\hat X,\ZZ) \lra H_3(X,\ZZ)\]
is surjective.\\

{\em Special role of the fourfold points}\\

So we see that we need to analyse the situation of a fourfold point
in more detail. We only have to consider fourfold points $p$ that are not
on triple lines, so at which four planes intersect in general position. 
Near $p$ the space $\RR^3$ 
is decomposed into $2^4=16$ cells, which are in one-to-one 
correspondence with the sign patterns 
\[ (\textup{sign}(L_1),\;\;\; \textup{sign}( L_2),\;\;\; \textup{sign} (L_3),\;\;\; \textup{sign}(L_4) )\]
where the linear forms define the planes meeting at $p$. Each cell
has an {\em opposite cell}, obtained by a reversing all signs.\\

\begin{center}
  \includegraphics[trim=6cm 10cm 4cm 10cm, clip, width=7cm]{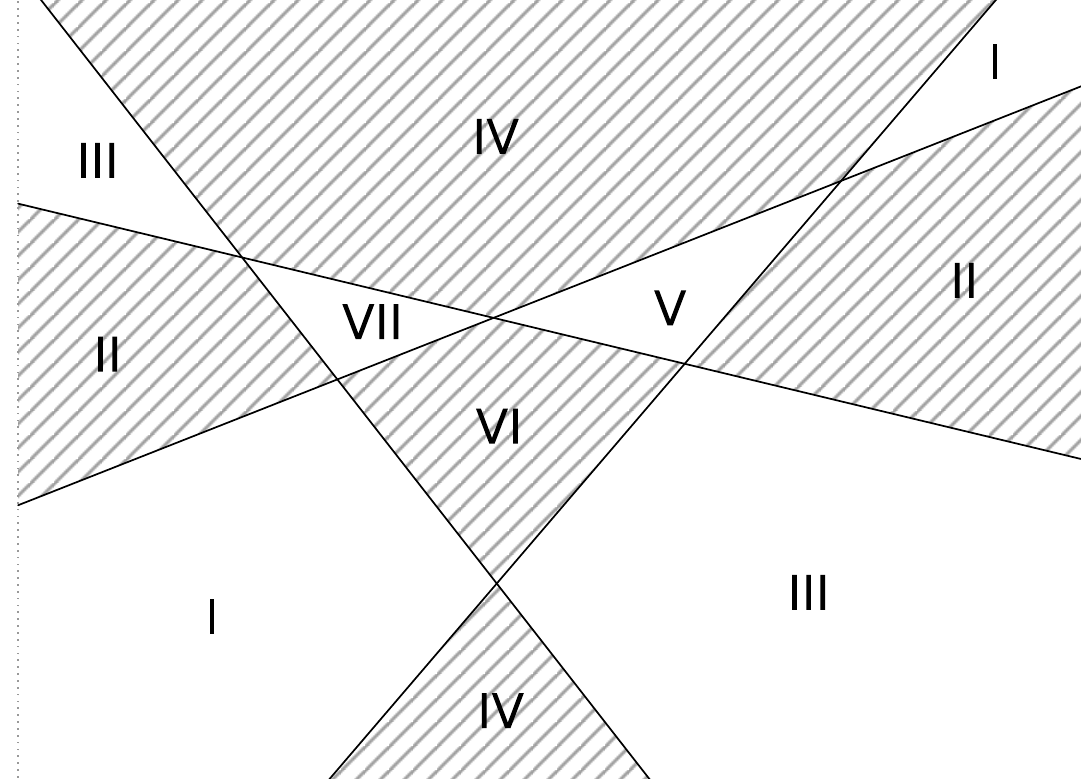}
\end{center}

If we blow--up the point $p$, the exceptional divisor is a copy of
$\PP[2]$, on which we find four lines in general position, 
corresponding to the four planes through $p$; these four lines decompose
the real projective plane in seven regions, that can be colored into
three 'black' and four 'white' regions.

\begin{center}
  \includegraphics[trim=5cm 10cm 3cm 10cm, clip, width=7cm]{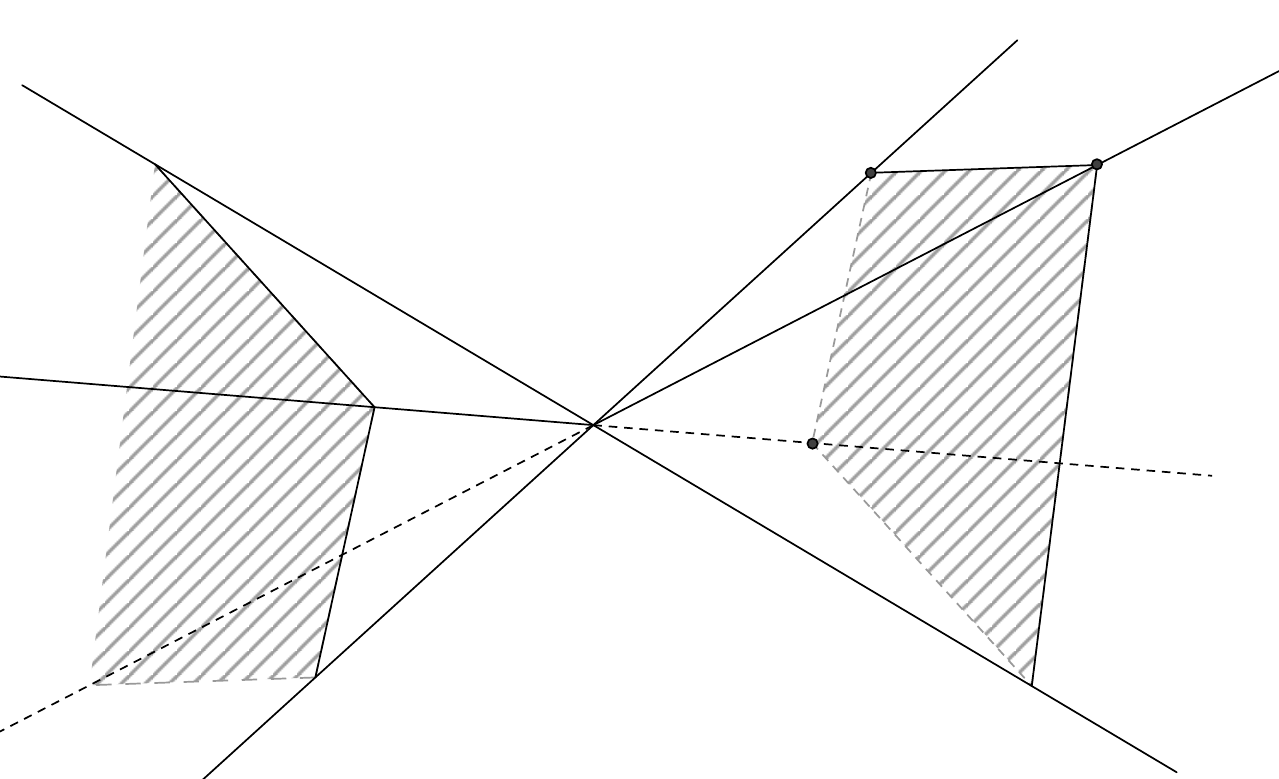}  
\end{center}

On the double cover we find an exceptional divisor $E$ that is a double
cover of this $\PP[2]$ ramified along these four lines, and each of the
eight regions $R$ determine a $2$-cycle ${\bf R}$ in $E$.

These regions are in one--to--one correspondence the pairs of opposite
cells. If $C$ is a cell corresponding to a region $R$, then the chain
${\bf C}$ on the blow--up, then the boundary of this chain is precisely
the $2$-cycle corresponding to $R$:
\[ \partial {\bf C} = \pm {\bf R}\]

What we learn from this is that we can cancel this boundary term of a cell
by adding to it the boundary term of the opposite cell!

Hence, we can define a group of {\em polyhedral cycles} $PC^{3}$ consisting
of elements
\[ \sum_C n_C {\bf C},\;\;\;n_C \in \ZZ\]
for which for all fourfold points $p$ one has:
\[ \sum_{ p \in \overline{C}} n_C =0\]


We can assume that the four planes we are
considering have equations 
\[x=0,\;\;\;y=0,z=0,\;\;\; x+y+z=0\] 
in appropriate coordinates. We will analyse what happens if we 
{\em smooth out} the fourfold point shifting the fourth plane to
\[x+y+z=\epsilon, \] 
resolve this and then specialise back to $\epsilon=0$.

Let us first blow--up the two disjoint lines 
\[ x=y=0,\;\;\; z=x+y+z-\epsilon=0. \]
In one of the affine charts the blow--up of $\mathbb P^{3}$ is given by
the equation 
\[x(y+1)-z(v-1)-\epsilon.\] 
The threefold is smooth unless $\epsilon=0$ when it acquires a node at
$x=0,y=-1,z=0,v=1$. Since the surface $x=0,z=0$ is  a Weil
divisor on the threefold which is not Cartier  (it is a
component of the exceptional locus of the blow--up) the node
admits a projective small resolution. Since the node does not lie on
the branch divisor it gives two nodes on the double cover.
Consequently. we get the partial resolution $\hat X$ from the end of the
section \ref{sec:doct}.

\section{Determination of $H_3(\hat X,\ZZ)$}

Denote by $X_{t}$ a smoothing of $\hat X$. By \cite{CvS}, the deformations of $\widetilde X$ correspond to deformations of the
arrangements of eight planes that preserve the incidences between the planes 
in $D$. The deformations of the arrangement that  preserves all the incidences 
except for the fourfold points correspond to {\em smoothings} of $\hat X$. 

By the work of {\sc J. Werner} \cite{Werner}, the nodal variety $\hat X$ 
is homotopy equivalent to its small resolution $\widetilde X$ with 3--cells glued  along the the exceptional lines (which topologically are $2$--spheres). 
The nodal variety $\hat X$ is also homotopy equivalent to its smoothing 
with $4$--cells glued along the vanishing $3$--cycles. As a consequence, one
arrives at the following equations relating topological invariants of
$X_t$, $\widetilde X$ and $\hat X$:
 
\begin{eqnarray*}
 && b_{4}(X_{t})+2p_{4}+b_{3}(\hat X)=b_{4}(\hat X)+b_{3}(X_{t})\\
 &&b_{2}(X_{t})=b_{2}(\hat X)\\
 && b_{3}(\widetilde X)+2p_{4}+b_{2}(\hat X)=b_{3}(\hat X)+b_{2}(\widetilde X)\\
 &&b_{4}(\widetilde X)=b_{4}(\hat X)
\end{eqnarray*}
and hence
\[b_{3}(\hat X)=b_{3}(\widetilde X)+b_{2}(X_{t})-b_{2}(\widetilde X)+2p_{4}=
b_{4}(\widetilde X)+b_{3}(X_{t})-b_{4}(X_{t})-2p_{4}
\]
where $p_{4}$ is the number of (smoothed) fourfold points in $D$ that
do not lie on a triple line. 

For the eleven rigid double octics from \cite{Meyer} we get
\[\def\arraystretch{1.3}
\begin{array}{|r||r|r|r|}
\hline
  \textrm{No.} & b_{3}(\hat X)&b_{3}(X_{t})&p_{4}^{0}\\
\hline
  1&3&4&1\\
  3&5&8&3\\
  19&6&10&4\\
  32&7&12&5\\
  69&7&12&5\\
  93&8&14&6\\
  238&11&20&12\\
  239&11&20&10\\
  240&11&20&10\\
  241&11&20&10\\
  245&11&20&9\\
\hline
\end{array}
\]

\section{An example}
In order to find two independent cycles, we draw projections of 
intersections of all arrangement planes onto the $(x,y)$--plane, 
and consider the equations of the planes not perpendicular to it 
as functions in $z$. The easiest case is the arrangement No 1. 
which has  a single $p_{0}^{4}$ point. We will go through some details
of this example . 

The equation of this arrangement is
\[xyzt(x+y)(y+z)(z+t)=0\]
and the only $p_{0}^{4}$ point is $(1,-1,1,-1)$.
The affine change of variables
\[t\longmapsto t-x,\]
maps this point to the plane at infinity. The arrangement is then given 
in affine coordinates by the equation
\[xyz(1-x)(x+y)(y+z)(-x+z+1)=0\]
while the $p_{0} ^{4}$ point is the point at infinity $(1,-1,1,0)$. 
The planes defined the first, second, fourth and fifth factor of the above product are perpendicular to the $(x,y)$--plane and intersect the plane in 
lines $x=0,y=0$,$x=1$,$x+y=0$. 
The planes in the arrangement that are not perpendicular to the $(x,y)$--plane 
can be seen as graphs over the $(x,y)$-plane and are given by
\begin{eqnarray*}
 &&z= f_{3}(x,y):=0\\
 &&z= f_{6}(x,y):=-y\\
 &&z= f_{7}(x,y):=x-1
\end{eqnarray*}
The projections of the lines of intersection of these planes are given by
\begin{eqnarray*}
  f_{3}=f_{6}&:&\qquad y=0\\
  f_{3}=f_{7}&:&\qquad x=1\\
  f_{6}=f_{7}&:&\qquad x+y=1
\end{eqnarray*}

In the $(x,y)$--plane we have two bounded domains
\begin{eqnarray*}
  I:&&\qquad x>0, y>0, x+y<1\\
  II:&&\qquad x<1, y<0, x+y>0
\end{eqnarray*}
\[\includegraphics[width=8cm]{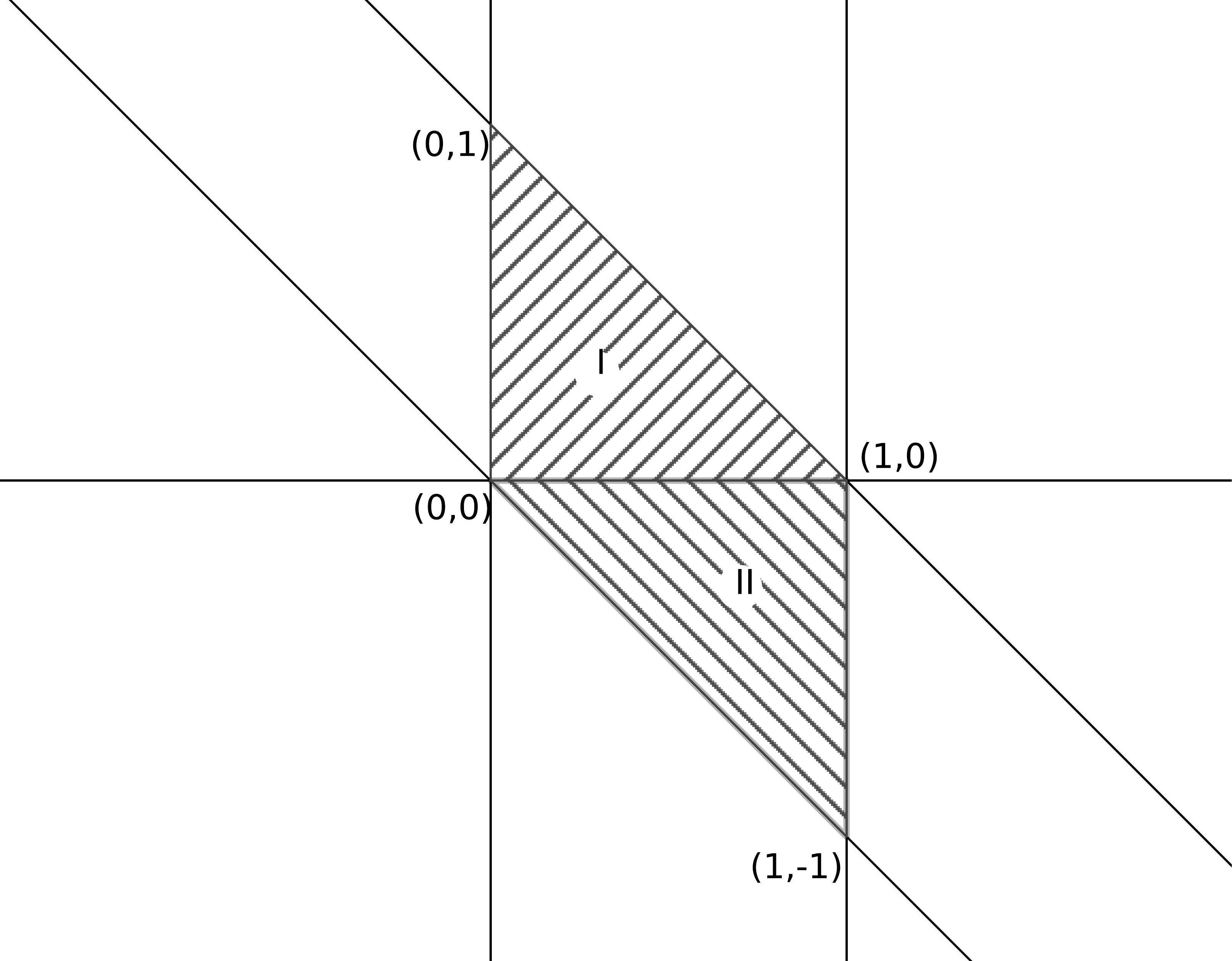}.\]

For points $(x,y)$ in these regions, the functions $f_{3}, f_{6}, f_{7}$ 
satisfy there the following inequalities
\begin{eqnarray*}
  I:&&\qquad f_{7}<f_{6}<f_{3}\\
  II:&&\qquad f_{7}<f_{3}<f_{6}
\end{eqnarray*}
Consequently, the domains lying over triangle I are given by
\begin{eqnarray*}
&&  x>0, y>0, x+y<1, z>x-1, z<-y\\
&&  x>0, y>0, x+y<1, z>-y,z<0
\end{eqnarray*}
As the ``right'' (horizontal) edge of the triangle $II$ is the projection
the intersection of planes no. 3 and 7, the only cycle lying over that
triangle is given by
\[x<1, y<0, x+y>0, z>x-1, z< 0,\]
the other domain 
\[x<1, y<0, x+y>0, z> 0, z<-y\]
is not bounded by arrangement planes, if we want to use it we would
have to add the unbounded domain ``across the edge'' 
\[x>1, x+y>0, x+y<1.\]

Instead we can choose a domain over triangle II
\[x>0, y<0, x+y>0, z>x-1, z<0.\]

\section{Period integrals}
When we are given a degree eight polynomial $F(x,y,z,t)$, then the
double octic $X \subset \PP[4](4,1,1,1,1)$ defined by the equation
\[ u^2-F(x,y,z,t)=0\]
comes with a preferred section $\omega \in \Gamma(X,\omega_X)$ of its
sheaf of dualising differentials. In the affine chart $t \neq 0 $ it can be 
written as 
\[\omega:= \frac{dxdydz}{u}=\frac{dxdydz}{\sqrt{F}}.\]

The period integrals of $X$ are thus of the form
\[\int_{\gamma} \omega=\int_{\gamma} \frac{dxdydz}{\sqrt{F}} \]
where $\gamma$ is a three-cycle in $X$.\\

If in particular $F$ defines a real arrangement of eight planes and
we have a bounded cell $C$ in $\mathbb{R^{3}}$ yielding a $3$--cycle 
$\bf{\widetilde{C}}$ in Calabi-Yau threefold $\widetilde X$, 
the period integral  
\[ \int_{\bf{\widetilde{C}}} \omega\]
is just equal to three-fold integral
\[2\iiint_{C}\frac{dx\;dy\;dz}{\sqrt{F}},\]

In the case considered  in previous section (arrangement No. 1), the two 
periods integral are given by
 
\begin{eqnarray*}
  &&\int _{0}^{1}\int_{0}^{1-x}\int_{x-1}^{-y}\frac1{\sqrt{xyz(1-x)(x+y)(y+z)(-x+z+1)}}dzdydx\\
  &&\int _{0}^{1}\int_{-x}^{0}\int_{x-1}^{0}\frac1{\sqrt{xyz(1-x)(x+y)(y+z)(-x+z+1)}}dzdydx.
\end{eqnarray*}

To compute such integrals numerically, we used {\sc Maple}.  
However, the function $F$ can have zeros of multiplicity $5$ at a vertex 
of a polyhedron of integration and thus the integrand is unbounded. 
As a result, a direct numerical integration usually does not yield a 
satisfactory precision in reasonable time. We used the following 
{\em simple trick} which allows us to get $12$ digits precision without 
much effort, which is sufficient for the our purposes. Using an affine
coordinate change, we reduce computations to the case of a integration over a cube 
$0 \le x,y,z \le 1$, with the function $F$ vanishing only for
$xyz=0$. Then substituting $(x,y,z) \mapsto (x^{k},y^{k},z^{k})$ in the triple 
integral transforms the integral to the integration of a {\em bounded} function.\\

Note that depending on the sign of the function $F$ in a given polyhedral cell $C$, we get either a real or a purely imaginary number. The computation time in the 
latter case can be reduced considerably by just using the function $-F$\, ! 
It should be noted that if we multiply $F$ by a constant factor $\lambda$, 
the corresponding  period integral changes by a factor $\sqrt{\lambda}$.
In particular, if we change the sign of $F$ the real and imaginary periods
are interchanged.\\

In the case of arrangement nr. $1$ everything works nicely, but in the other 
cases the picture of the decomposition of $\mathbb P^{3}$ becomes much more 
complicated and more generic fourfold points to take into account. 
We wrote a simple {\sc Maple} code to produce a linear--cylindric 
decomposition and form cycles from the polyhedral cells. Then
we used several changes of variables moving each of the planes of the
arrangement to infinity, which allowed us to compute the integrals for
most of the cycles. In all cases ratios of any two real and any two complex
integrals were rational numbers (with numerator and denominator
$\le 6$). Below is a table that summarises all different period integrals
that appeared in our calculations.\\

It should be kept in mind that all period integrals get multiplied by
a common factor if we change the polynomial $F$ defining the
arrangement.
For these calculations we used the equations $F$ as listed in
\cite{Meyer} scaled by $\lambda$ from the last column of table at the
end of section~\ref{sec:doct}.\\

{\scriptsize
\begin{tabular}{r||l|l}
  Arr. No.& Real integrals& Imaginary integrals\\\hline
  1&55.9805041334, 111.961008267 &69.3694986501i\\
  3&80.3028893419, 160.60577868&41.4134587444i, 82.8269174889i \\
  &&124.240376233i,  289.89421121i\\
  19&72.1085316451, 144.217063291&72.1085316451i, 144.217063291i\\
  &216.325594935&216.325594935i\\
  32&55.9805041335,   111.961008267&34.6847493250i, 69.3694986501i\\
  &&138.738997300i, 208.1084959i\\
  69&55.9805041335, 111.9610083&34.6847493252i, 138.738997300i\\
  &223.922016533&277.4779945i\\
  93&55.9805041334&17.3423746625i, 69.3694986502i\\
  &&138.738997300i\\
  238&55.9805041334, 111.961008267 &34.6847493250i\\
  239&48.5252148713, 145.575644614&35.2275632784i, 105.682689835i\\
  240&43.7468074540, 131.240422363&28.8234453872i, 57.6468907743i\\
  241&223.922016533&69.3694986503i\\
  245&21.8734037270, 87.4936149079&28.8234453872i, 115.293781548i\\
  &131.240422362&
\end{tabular}
}
\smallskip

For each arrangement, the computed period integrals generate 
a lattice in $\CC$, which in turn defines an elliptic curve.
This lattice might be a proper sublattice of $H_3(\widetilde X,\ZZ)$,
but in any case it defines a elliptic curve that is 
 isogeneous with the intermediate Jacobian $J^2(\widetilde X)$ 
of the corresponding Calabi--Yau threefold. In the following table we 
list lattice generators, $j$--invariant and coefficients of the 
classical Weierstrass equation 
\[ y^2=4 x^3-g_2 x-g_3 \]
of the elliptic curve, which are easily computed 
numerically via
\[g_2=60 \sum_{0 \neq m \in \Lambda}\frac{1}{m^4}, \;\;\; 
g_6= 140 \sum_{0 \neq m \in \Lambda}\frac{1}{m^6}, \;\;\; 
\]
This is a standard functionality in {\sc MAGMA}. 
\smallskip

{\scriptsize
\begin{tabular}{r||r|r|r|r}
  Arr. No&$\tau/i$&$j(\tau)$&$g_{2}$&$g_{3}$\\\hline
1&1.23917245341& 3236.13720434&142.879810750&224.378572683\\
3&0.515715674539& 196267.167917&1838.35630102&-15102.274126\\
19&1&1728&189.072720130&0\\
32&0.619586226703&26112.0318779&889.658497527&-4934.98162416\\
69&0.619586226703&26112.0318779&889.658497527&-4934.98162416\\
93&0.309793113352&643142260.966&14101.0467615&-322251.215146\\
238&0.619586226704&26112.0318791&889.658497527&-4934.98162416\\
239&0.725964086338&6517.46790207&487.190579154&-1774.06947556\\
240&0.658869688205&14612.0507801& 701.139041736& -3355.01890381\\
241&0.309793113354&643142256.756&14101.0467615&-322251.215146\\
245&1.31773937641&4737.95402281& 137.802991416& 248.136467781
\end{tabular}
}
\smallskip

\section{Comparison with modular periods}
Recall that for an Hecke-eigenform $f \in S_k(\Gamma_0(N))$
with $q$-expansion  
\[ f=\sum_{n=1}a_n q^n\]
the $L$-function is defined by the series:
\[ L(f,s)=\sum_{n=1}^\infty \frac{a_n}{n^s} . \]
It converges for $Re(s) >1+k/2$ and the completed $L$-function
\[\Lambda(f,s):=(\sqrt{N}/2\pi)^s \Gamma(s)L(f,s)\]
satisfies the functional equation
\[ \Lambda(f,s)=w i^k \Lambda(f,k-s)\]
where $w$ is the sign of $f$ under the Atkin-Lehner involution.
We note that 
\[\Lambda(f,s) =(\sqrt{N})^s \int_0^{\infty}f(it) t^s \frac{dt}{t}\]

In our case we $k=4$ and $w=1$, so that the functional
equation just reads
\[ \Lambda(f,s)=\Lambda(f,4-s)\]
This means in particular
\[\Lambda(f,1)=\Lambda(f,3)\]
from which we get the equality
\[L(f,3) = \frac{(2\pi)^2}{N} \frac{\Gamma(1)}{\Gamma(3)} L(f,1)=\frac{2\pi^2}{N}L(f,1)\]
Furthermore, we see from the functional equation that $L(f,k)=0$ for $k=0,-1,-2,-3,-4,\ldots$.

By direct point counting (and correcting for the singularities of course),
C. Meyer was able to determine cusp-forms 
\[f=\sum_{n=1}^{\infty} a_n q^n \in S_4(\Gamma_0(N))\]
such that
\[ a_p =Tr(Fr_p:H^3(X) \lra H^3(X))\]
In other words, one has equality of $L$-functions
\[L(H^3(X),s)=L(f,s) \]
The result is summarised in the following table (we multiplied the
equation of the octic arrangement by the factor $\lambda$ to obtain
modular form of minimal level).
\smallskip

{\scriptsize
\[
\begin{array}{|r|l|l|}
 \hline
\textup{Form}&\textup{q-expansion}&\textup{Arrangements}\\
\hline
  6/1&q - 2q^2 - 3q^3 + 4q^4 + 6q^5 + 6q^6 - 16q^7 + O(q^{8})    &240, 245\\
  8/1&q - 4q^3 - 2q^5 + 24q^7 - 11q^9 - 44q^{11} + O(q^{12})   &1, 32, 69, 93, 238, 241\\
  12/1&q + 3q^3 - 18q^5 + 8q^7 + 9q^9 + 36q^{11} + O(q^{12}) &239\\
  32/1&q + 22q^5 - 27q^9 + O(q^{12})   &19\\
  32/2&q + 8q^3 - 10q^5 + 16q^7 + 37q^9 - 40q^{11} + O(q^{12})&3\\
\hline
\end{array}
\]}

\noindent
It is a remarkable fact that only five different modular forms appear.

In cases where two varieties $X, X'$ have the same modular form, one expects 
there exists a correspondence $\phi$ between $X$ and $X'$ that explains it.

It is gratifying to see that the numerical evaluation of the period integrals
lead to the very same grouping of our examples. 


Here we summarise the calculations of the critical $L$-values  
\[\scriptsize
\begin{array}{|c||c|c|}
\hline
f&L(f,1)&L(f,2)\\
\hline
6/1&0.22162391559067350824671004425&0.50971042336159397988737819140\\
\hline
8/1&0.35450068373096471876555989149&0.69003116312339752511910542021\\
\hline
12/1&0.61457902590673022954002802969&0.93444013814191444281042898230\\
\hline
32/1&1.82653044425089816105284840591&1.43455365630418076432004680798\\
\hline
32/2&2.03409594950627923591429024672&1.64778916742512594127684239683\\
\hline
\end{array}
\]

\smallskip

If we express the real and imaginary periods of the double octics we
get, at least at the numerical level, nice proportionalities with
\[\pi L(f,2),\;\;\; \pi^2 L(f,1)\]
for the corresponding modular form.

\[
\begin{array}{|c|c|c|}
\hline
\multicolumn{3}{|c|}{\text{Form 6/1}}\\
\hline
\hline
240& 43.7468074540\ldots=20\pi^{2} L(f,1)&28.8234453871\ldots=18\pi L(f,2)\\
245& 21.8734037270\ldots=10\pi^{2} L(f,1)&28.8234453871\ldots=18\pi L(f,2)\\
\hline
\hline
\multicolumn{3}{|c|}{\text{Form 8/1}}\\
\hline
\hline
 1&55.9805041334\ldots=16\pi^2 L(f,1)&69.3694986501\ldots=32\pi L(f,2)\\
32&55.9805041334\ldots=16\pi^2 L(f,1)&34.6847493250\ldots=16\pi L(f,2)\\
69&55.9805041334\ldots=16\pi^2 L(f,1)&34.6847493250\ldots=16\pi L(f,2)\\
93&55.9805041334\ldots=16\pi^2 L(f,1)&17.3423746625\ldots=8\pi L(f,2)\;\;\\
238&55.9805041334\ldots=16\pi^2L(f,1)&34.6847493250\ldots=16\pi L(f,2)\\
\hline
\hline
\multicolumn{3}{|c|}{\text{Form 12/1}}\\
\hline
\hline
239&48.5252148713\ldots=8 \pi^2 L(f,1)&35.2275632785\ldots=12\pi L(f,2)\\
\hline
\hline
\multicolumn{3}{|c|}{\text{Form 32/1}}\\
\hline
\hline
19&72.1085316452\ldots=4\pi^{2} L(f,1)&72.1085316452\ldots=16\pi L(f,2)\\
\hline
\hline
\multicolumn{3}{|c|}{\text{Form 32/2}}\\
\hline
\hline
3&80.3028893419\ldots=4\pi^2 L(f,1)&41.4134587443\ldots=8\pi L(f,2)\\
\hline
\end{array}
\]

\smallskip

\smallskip

\section{Outlook}
The above calculations show that it is possible to verify numerically
the relation between the periods of a rigid Calabi-Yau and the 
corresponding $L$-values of the attached modular form. However, one
would like to push these calculations to a further level. 
One important problem that was left untouched by our calculations is
the complete determination of the group $3$-cycles in terms
of polyhedral cycles. We identified some polyhedral cycles, but there is 
no guaranty that these generate the whole third homology group 
$H_3(\tilde X,\ZZ)$. It follows from Poincar\'e-duality, that this group is generated by any two cycles with intersection $\pm1$. This leads to the following
question to determine the inteeresection number $\langle \delta, \gamma \rangle$
between two polyhedral cycles $\delta, \gamma$ in purely combinatorial
terms of the cell appearing in $\delta$ and $\gamma$, and their mutual
position inside the arrangement. 


Apart from the eleven double octics with rational coefficients there
are two other with coefficients in  $\QQ(\sqrt{-3})$ and one in
$\QQ(\sqrt{5})$ (cf. \cite{CvS2,CvS3}). However, in these cases
presented method did not yield any reasonable approximation of the
period integrals.\\

\textbf {Acknowledgments:} Part of this research was done during the 
stay of the first named author  as a guestprofessor at the \emph{Schwerpunkt
  Polen} of the Johannes  Gutenberg--Universit\"at in Mainz. He would
like to thank the department for the hospitality and excellent working
conditions. This research was supported in part by PL--Grid
Infrastructure.  The first named author was partially supported by NCN
grant no. N N201 608040.

\end{document}